\documentclass[12pt,twoside]{amsart}

\newtheorem{Thm}{Theorem}
\newtheorem{Cor}[Thm]{Corollary}
\newtheorem{Lemma}[Thm]{Lemma}
\newtheorem{Prop}[Thm]{Proposition}
\newtheorem*{Prop*}{Proposition}

\theoremstyle{definition}

\title{K-homology of the CAR algebra}
\author[T.~Hadfield]{Tom Hadfield}
\address{Department of Mathematics, University of California, Berkeley, CA 94720}
\email{hadfield@math.berkeley.edu}
\date{\today}

\begin{document}

\begin{abstract}
 We compare the K-theory and K-homology of the well known CAR algebra. 
 While the $K_0$ group is infinitely generated, 
we show that it pairs identically to zero with even Fredholm modules 
over the algebra. 
We show further that in fact the even K-homology is trivial, 
while the odd K-homology group is very large. 
 These results give some insight into the K-homology of AF-algebras. 
\end{abstract}

\maketitle

\section{The CAR algebra}
 
We begin by recalling the definition of the CAR algebra. 
For each positive integer $n$, denote by $A_n$ the full matrix algebra $M_{2^n}({\bf C})$, and  
let $\phi_{n+1,n} : A_n \rightarrow A_{n+1}$ be the *-homomorphism 
\begin{equation}
 [B_n] \mapsto 
\left[
\begin{array}{cc}
B_n & 0 \cr
0 & B_n \cr
\end{array}
\right].
\end{equation}
 For each $k \geq 1$, we define $\phi_{n+k,n} : A_n \rightarrow A_{n+k} $
by $\phi_{n+k,n} = \phi_{n+k,n+k-1} \circ \ldots \circ \phi_{n+1,n}$. 
 Let $A = {lim_{\rightarrow}}{A_n}$ be the CAR algebra, and let 
 $\phi_n : A_n \rightarrow A$ be the inclusion maps. 
 For each $k \geq 0$ we have 
$\phi_n = \phi_{n+k} \circ {\phi_{n+k,n}}$. 

 The K-theory of the CAR algebra is well known \cite{davidson} :

\begin{Prop}
${K_0}(A) \cong {\bf Z}[ {\frac{1}{2}}]$, ${K_1}(A) =0$. 
\end{Prop}

For each $n$, we have 
${K_0}(A_n) \cong {\bf Z}$, generated by $[p_n]$, where $p_n \in A_n$ is a rank one projection, and  ${K_1}(A_n) =0$.

\section{K-homology}

 Since the $A_n$ are full matrix algebras their K-homology is well known. 
We have 
${KK^0}(A_n, {\bf C})$ $\cong {\bf Z}$ and ${KK^1}(A_n, {\bf C}) \cong 0$. 

\begin{Lemma}
 The generator of the even K-homology of $A_n$ is the canonical even Fredholm module
\begin{equation}
{\bf z}_0 ^{(n)} = ( {\bf C}^{2^n} \oplus {\bf C}^{2^n}, \pi_n = id \oplus 0 , 
F= \left(
\begin{array}{cc}
0 & 1 \cr
1 & 0 \cr
\end{array}
\right),
\gamma =
\left(
\begin{array}{cc}
1 & 0 \cr
0 & -1 \cr
\end{array}
\right)
)
\end{equation}
\end{Lemma}
 
\begin{proof} It is easy to check that ${\bf z}_0 ^{(n)}$ is a Fredholm module. 
 It is immediate that $F = {F^{*}}$, $F^2 =1$, and the commutators $[F, \pi(a)]$ are all trivially compact, since we are working with finite dimensional Hilbert spaces. 
 To show nontriviality, we show that the Chern character of ${\bf z}_0 ^{(n)}$ pairs to 1 with the generator $[p_n]$ of $K_0 (A_n)$. 
Recall that
\begin{equation}
 < {ch_{*}}( {\bf z}_0^{(n)} ), [p_n] > = 
{{\lim}_{k \rightarrow \infty} } {(k!)}^{-1} \psi_{2k}  (p_n ,..., p_n).
\end{equation}
 Here  
${ch}_{*} : {KK^0}(A_n ,{\bf C}) \rightarrow {HC^{even}}(A_n )$, 
is the even Chern character as defined in \cite{connes94}, p295, mapping the even K-homology of $A_n$ into even periodic cyclic cohomology,
 $<.,.>$ denotes the pairing between K-theory and periodic cyclic cohomology defined in \cite{connes94}, p224, and 
 $\psi_{2k}$ is the cyclic $2k$-cocycle given by 
\begin{equation}
 \psi_{2k} ( a_0, a_1, ..., a_{2k}) = 
(-1)^{k(2k-1)} \Gamma (k+1) 
 Tr( \gamma \pi(a_0) [ F, \pi(a_1)] ... [F, \pi(a_{2k})]) 
\end{equation}
Since $\Gamma(k+1) = k!$ we have 
\begin{equation}
< {ch_{*}}( {\bf z}_0^{(n)} ), [p_n] > = 
{{\lim}_{k \rightarrow \infty} }
 {(-1)}^k
 Tr( \gamma \pi_n (p_n) {[F, \pi_n (p_n)]}^{2k})
\end{equation} 
Now,  
$[F, \pi_n ( p_n ) ] = 
 \left(
\begin{array}{cc}
0 & - p_n \cr
p_n & 0 \cr
\end{array}
\right)$,
 hence
$ \gamma \pi_n ( p_n) {[F, \pi_n (p_n)]}^{2k} =
(-1)^k 
\left(
\begin{array}{cc}
p_n & 0 \cr
0 & 0 \cr
\end{array}
\right).$
 Therefore 
\begin{equation}
 < {ch_{*}}( {\bf z}_0^{(n)} ), [p_n] > = 
{{\lim}_{k \rightarrow \infty} } {(-1)}^k 
Tr( 
(-1)^k 
\left(
\begin{array}{cc}
p_n & 0 \cr
0 & 0 \cr
\end{array}
\right)) = 1,
\end{equation}
as claimed. 
\end{proof}

\begin{Lemma}
${\phi_{n+1,n}^{*}} ({{\bf z}_0 ^{(n+1)}}) = 2{{\bf z}_0 ^{(n)}}$.
\end{Lemma}

\begin{proof}
 By definition, 
\begin{equation}
{\phi_{n+1,n}^{*}} ({{\bf z}_0 ^{(n+1)}}) =
( {\bf C}^{2^{n+1}} \oplus {\bf C}^{2^{n+1}},  id \circ \phi_{n+1,n} \oplus 0 , 
F= \left(
\begin{array}{cc}
0 & 1 \cr
1 & 0 \cr
\end{array}
\right),
\gamma =
\left(
\begin{array}{cc}
1 & 0 \cr
0 & -1 \cr
\end{array}
\right)
)
\end{equation}
\begin{equation} 
=
( ({\bf C}^{2^n} \oplus {\bf C}^{2^n}) \oplus 
({\bf C}^{2^n} \oplus {\bf C}^{2^n})
, (id \oplus id) \oplus 0 , 
F= \left(
\begin{array}{cc}
0 & 1 \cr
1 & 0 \cr
\end{array}
\right),
\gamma =
\left(
\begin{array}{cc}
1 & 0 \cr
0 & -1 \cr
\end{array}
\right)
)
\cong 
{{\bf z}_0 ^{(n)}} + {{\bf z}_0 ^{(n)}}
\end{equation}
\end{proof}

Hence the tower of abelian groups 
\begin{equation}
\ldots \rightarrow {KK^0}( A_n , {\bf C}) \rightarrow \ldots \rightarrow 
{KK^0}( A_1 , {\bf C}) \rightarrow
{KK^0}( A_0 , {\bf C})
\end{equation} 
is just 
$\ldots \rightarrow {\bf Z} \rightarrow \ldots \rightarrow {\bf Z} \rightarrow {\bf Z} $
 with each connecting map being given by $m \mapsto 2m$. 
It follows that :

\begin{Lemma}
\label{inverselimit}
The inverse limit of the even K-homology groups is trivial, 
 ${\lim_{\leftarrow}} {KK^0}( A_n, {\bf C}) \cong 0.$  
\end{Lemma}
\begin{proof}
 Recall that, for a tower of abelian groups 
\begin{equation}
 \ldots \rightarrow G_n \rightarrow^{f_n} G_{n-1} \rightarrow^{f_{n-1}} \ldots \rightarrow^{f_0} G_0
\end{equation}
 the inverse limit ${\lim_{\leftarrow}}{G_n}$ is isomorphic to the abelian group consisting of all sequences ${\{ g_n \}}_{n \geq 0}$, with $g_n \in G_n$ for each $n$, such that $g_{n-1} = {f_n}(g_n)$. 
 In our situation, an element of the inverse limit will be a sequence $\{ a_n \}$, $a_n \in {\bf Z}$, such that $a_n =2 a_{n+1}$ for each $n$. Hence we must have all $a_n =0$. 
 So the inverse limit is the trivial group. 
\end{proof}

Since the K-theory of the CAR algebra $A$ is not free abelian, 
Rosenberg and Schochet's 
universal coefficient theorem \cite{rs}, \cite{blackadar} p234 does not give us the K-homology for free. 
However :

\begin{Prop} 
Given any ${\bf z} \in {KK^0}(A,{\bf C})$, we have 
${\phi_n^{*}}({\bf z}) = {\bf 0} \in {KK^0}( A_n, {\bf C})$.
\end{Prop}
\begin{proof}
Given $[x] \in {K_0}(A)$, since $A$ is AF, there exists $N$ so that for each $n \geq N$ there is  
$[x_n] \in {K_0}(A_n)$, 
such that 
${\phi_n}_{*}[x_n] = [x]$. 
Consider 
${\phi_n^{*}}({\bf z}) \in {KK^0}(A_n, {\bf C})$. 
 We must have 
${\phi_n^{*}}({\bf z}) = {k_n}  {\bf z}_0^{(n)}$ 
for some 
${k_n} \in {\bf Z}$. 
 Suppose 
${k_n} \neq 0$. 
Then we can choose 
$p \in {\bf Z}_{>0}$ so that $2^p$ does not divide $k_n$. 
 Now, we have 
$\phi_n = \phi_{n+p} \circ \phi_{n+p,n}$. 
 So 
${\phi_n^{*}}( {\bf z}) = {\phi_{n+p,n}^{*}}( {\phi_{n+p}^{*}} ( {\bf z}))$, 
and  
${\phi_{n+p}^{*}} ( {\bf z}) = {k_{n+p}} {{\bf z}_0^{(n+p)}}$. 
Hence 
${\phi_{n}^{*}}( {\bf z}) = {2^p} {\phi_{n+p}^{*}} ( {\bf z})$, 
so
${k_n} = 2^p k_{n+p}$. 
But $k_n$ is by assumption a nonzero integer, not divisible by  $2^p$.
 So this is a contradiction. 
So we must have $k_n = 0$ for all $n$, thus 
${\phi_n^{*}}({\bf z}) = {\bf 0} \in {KK^0}( A_n, {\bf C})$.  
\end{proof}

\begin{Cor}
Given any 
${\bf z} \in {KK^0}(A,{\bf C})$, 
and any 
$[x] \in {K_0}(A)$, 
we have 
$<{ch_{*}}( {\bf z}), [x]>=0$. 
\end{Cor}

\begin{proof}
 Given $[x] \in {K_0}(A)$, choose $n$ and $[x_n] \in {K_0}(A_n)$ such that 
${\phi_n}_{*}[x_n] = [x]$.  Then \\
$ < {ch_{*}}( {\bf z}), [x]> =
< {ch_{*}}( {\bf z}), {\phi_n}_{*}[x_n]> = 
< {ch_{*}}( {\phi_n^{*}} ({\bf z})), [x_n]> = 
< {ch_{*}}( {\bf 0}), [x_n]> =0$.
\end{proof}

We would like to deduce from this that ${KK^0}(A,{\bf C})=0$, but we do not know whether the pairing $<.,.>$ of K-theory and K-homology is faithful. 
However, we can use the following special case of a much more general result of Rosenberg and Schochet \cite{rs} :

\begin{Prop}
\label{AFkhom}
Suppose that $A = {lim_{\rightarrow}} A_n$ is an AF-algebra. Then  the following sequences on K-homology are exact:
\begin{equation*}
0 \rightarrow {lim^1}_{\leftarrow}  {KK^1}( A_n, {\bf C}) 
\rightarrow {KK^0}(A,{\bf C}) 
\rightarrow {lim_{\leftarrow}} {KK^0}( A_n , {\bf C}) 
\rightarrow 0
\end{equation*}
\begin{equation*}
0 \rightarrow 
{lim^1}_{\leftarrow}  {KK^0}( A_n, {\bf C}) 
\rightarrow {KK^1}(A,{\bf C}) 
\rightarrow {lim_{\leftarrow}} {KK^1}( A_n , {\bf C}) 
\rightarrow 0
\end{equation*}
\end{Prop}

The left hand term is Milnors' ${lim^1}_{\leftarrow}$, and the right hand term is the inverse limit of the  K-homology groups.  
We showed (Lemma \ref{inverselimit}) that 
 ${lim_{\leftarrow}} {KK^0}( A_n , {\bf C}) \cong 0$. 
 Furthermore, since ${KK^1}( A_n , {\bf C}) \cong 0$ for each $n$, we have 
 ${lim_{\leftarrow}} {KK^1}( A_n , {\bf C}) \cong 0$ also. 

 We have ${KK^0}(A_n, {\bf C}) \rightarrow {KK^0}(A_{n-1}, {\bf C})$ corresponding to ${\bf Z} \rightarrow {\bf Z}$, $m \mapsto 2m$.
 Equivalently, this is 
\begin{equation}
{\bf Z} \hookleftarrow 2 {\bf Z} \hookleftarrow ...
\hookleftarrow 2^n {\bf Z} \hookleftarrow ...
\end{equation}
 with all the maps just being inclusions. 
 This is exactly the situation of \cite{weibel}, Example 3.5.5, p82.  
 We therefore  have 
 ${lim^1}_{\leftarrow}{KK^0}( A_n, {\bf C})$ $={\hat{\bf Z}}_2 / {\bf Z}$, 
 where ${\hat{\bf Z}}_2$ is the uncountable (additive) group of 2-adic integers . 
 Hence for the CAR algebra $A$, we have 
${KK^1}(A,{\bf C}) \cong$ 
${{lim^1}_{\leftarrow}} {KK^0}( A_n , {\bf C})=$
$={\hat{\bf Z}}_2 / {\bf Z}$. 
 
The maps 
${KK^1}(A_n, {\bf C}) \rightarrow {KK^1}(A_{n-1}, {\bf C})$ 
 are all just the trivial surjective maps $0 \rightarrow 0$, 
hence from  \cite{weibel}, p80, 
${lim^1}_{\leftarrow}  {KK^1}( A_n, {\bf C}) = 0$ also.
 So ${KK^0}(A,{\bf C}) \cong$ 
${lim_{\leftarrow}} {KK^1}( A_n , {\bf C})=0.$ 
  
 Therefore we have shown that:

\begin{Prop}
\label{khomCAR}
Let $A$ be the CAR algebra. Then ${KK^0}(A,{\bf C}) \cong 0$, while\\  
 ${KK^1}(A,{\bf C})$ $\cong{\hat{\bf Z}}_2 / {\bf Z}$.
\end{Prop}

We note that these results are in agreement with \cite{blackadar}, Exercises 16.4.7, p141, and 23.15.2, p245.

\section{Acknowledgements}
 I would like to thank my advisor Professor Marc Rieffel for his advice and support throughout my time in Berkeley. I am extremely grateful for his help.

\bibliographystyle{amsalpha}

\begin{thebibliography}{1}

\bibitem[Bla98]{blackadar} {\sc b. blackadar} \emph{K-theory for operator algebras} : MSRI publications, 5. Cambridge University Press (1998).

\bibitem[Co94]{connes94} {\sc a. connes} \emph{Noncommutative geometry} : Academic Press (1994).

\bibitem[Da96]{davidson} {\sc k. davidson } \emph{C*-algebras by example} : Fields Institute Monographs {\bf 6} , Amer. Math. Soc. , Providence, RI (1996).

\bibitem[RS87]{rs} {\sc j. rosenberg c. schochet} \emph{The Kunneth theorem and the universal coefficient theorem for Kasparov's generalised K-functor} : Duke Math. J. {\bf 55}:2 , 431-474 (1987).

\bibitem[Wei94]{weibel} {\sc c. weibel} \emph{An introduction to homological algebra} : Cambridge University Press (1994).

\end{thebibliography}

\end{document}